\documentclass[12pt,letterpaper]{article}
\usepackage{amsfonts}
\usepackage{amssymb}
\usepackage{enumerate}

\newtheorem{ttt}{Theorem}[section]
\newtheorem{llll}[ttt]{Lemma}
\newtheorem{ccc}[ttt]{Claim}
\newtheorem{eee}[ttt]{Example}
\newtheorem{sss}[ttt]{Statement}
\newtheorem{ddd}[ttt]{Definition}
\newtheorem{qqq}[ttt]{Question}
\newtheorem{cccc}[ttt]{Corollary}

\newcommand{\bt}{\begin{ttt}}
\newcommand{\bl}{\begin{llll}}
\newcommand{\bc}{\begin{ccc}}
\newcommand{\bex}{\begin{eee}}
\newcommand{\bs}{\begin{sss}}
\newcommand{\bd}{\begin{ddd} \upshape}
\newcommand{\bq}{\begin{qqq}}
\newcommand{\bcor}{\begin{cccc}}

\newcommand{\et}{\end{ttt}}
\newcommand{\el}{\end{llll}}
\newcommand{\ec}{\end{ccc}}
\newcommand{\eex}{\end{eee}}
\newcommand{\es}{\end{sss}}
\newcommand{\ed}{\end{ddd}}
\newcommand{\eq}{\end{qqq}}
\newcommand{\ecor}{\end{cccc}}

\newcommand{\RR}{\mathbb{R}}

\newcommand{\fc}{\mathfrak{c}}

\newcommand{\al}{\alpha}

\newcommand{\iF}{\mathcal{F}}

\newcommand{\iP}{\mathcal{P}}

\newcommand{\sms}{separable metric space}

\newcommand{\lh}{\mathrm{length}}  
\newcommand{\cf}{\mathrm{cf}}  

\newcommand\cat{^{\mathord{\frown}}}

\newcommand{\CH}{\mathit{CH}}

\newcommand{\MA}{\mathit{MA}}
\newcommand{\ZFC}{\mathit{ZFC}}

\newcommand\rest{\mathord{\upharpoonright}}   

\newcommand\eop{$\ \ {\vcenter
   {\hrule
   \hbox{\vrule height 9pt \kern 9pt \vrule height 9pt}
   \hrule}}$\vskip 2 pt}

\newenvironment{proof}{{\bf Proof.}}{\hfill\eop\medskip}

\newcommand{\bp}{\begin{proof} }
\newcommand{\ep}{\end{proof} }

\begin{document}

\title{Transfinite Sequences of Continuous and Baire Class 1 Functions}

\author{M\'arton Elekes
\  and Kenneth Kunen\thanks
{Partially supported by NSF Grant DMS-0097881.}}

\maketitle 

\begin{abstract}
The set of continuous or Baire class 1 functions defined on a metric space
$X$ is
endowed with the natural pointwise partial order. We investigate how the
possible lengths of well-ordered monotone sequences (with respect to
this order) depend on the space $X$. 
\end{abstract}

\insert\footins{\footnotesize{

}}

\section*{Introduction}
Any set $\iF$ of real valued functions defined on an arbitrary set $X$ is
partially ordered by the pointwise order; that is, $f\leq g$ iff
$f(x)\leq g(x)$ for all $x\in X$.  Then,
$f<g$ iff $f \le g$ and $g \not\le f$; equivalently, $f(x)\leq 
g(x)$ for all $x\in X$ and $f(x) < g(x)$ for at least one $x\in X$. Our
aim will be to investigate the possible lengths of the increasing or
decreasing well-ordered sequences of functions in $\iF$ with respect to
this order. A classical theorem
(see Kuratowski \cite{Ku}, \S 24.III, Theorem $2'$)
asserts that if $\iF$ is the
set of Baire class 1 functions
(that is, pointwise limits of continuous functions)
defined on a Polish space $X$ (that
is, a complete separable metric space), then there exists a monotone
sequence of length $\xi$ in $\iF$ iff $\xi<\omega_1$.
 P. Komj\'ath \cite{Ko} proved that the corresponding
question concerning Baire class $\al$ functions for $2\leq\al<\omega_1$ is
independent of $\ZFC$. 

In the present paper we investigate what happens if we
replace the Polish space $X$ by an arbitrary metric space.

Section \ref{sec-cont} considers chains of continuous functions.
We show that for any metric space $X$, there exists a chain
in $C(X,\RR)$ of order type $\xi$ iff  $|\xi|\leq d(X)$.
Here, $|A|$ denotes
the cardinality of the set $A$, while $d(X)$ denotes the density of the space
$X$, that is
\[
d(X) = \max(\min\{|D|: D\subseteq X \ \&\ \overline{D} = X\}, \  \omega ) \ \ .
\]
In particular, for separable $X$, every well-ordered chain has
countable length, just as for Polish spaces.

Section \ref{sec-Baire} considers chains of Baire class 1 functions
on separable metric spaces.  Here, the situation is entirely
different from the case of Polish spaces, since
on some separable metric spaces, there are well-ordered
chains of every order type less than $\omega_2$.
Furthermore, the existence of chains of type $\omega_2$ and
longer is independent of $\ZFC + \neg \CH$.
Under $\MA$, there are chains of all types less than $\fc^+$,
whereas in the Cohen model, all chains have type less than $\omega_2$.

We note here that instead of examining well-ordered sequences, which is a
classical problem, we could try to characterize all the possible order
types of linearly ordered subsets of the partially ordered set
$\iF$. This problem was posed by M. Laczkovich,
and is considered in detail in \cite{El}.

\section{Sequences of Continuous Functions}
\label{sec-cont}

\bl
For any topological space $X$:
If there is a well-ordered sequence of length $\xi$ in $C(X,\RR)$,
then  $\xi<d(X)^+$.
\el
\bp
Let $\{ f_\al :\al<\xi\}$ be
an increasing sequence in $C(X,\RR)$, and let
$D\subseteq X$ be a dense subset of $X$ such that
$d(X) = \max(|D|,\omega)$.
By continuity, the $f_\alpha \rest D$ are all distinct; so, for each
$\al<\xi$,  choose a $d_\al\in D$ such that
$f_\al(d_\al) < f_{\al+1}(d_\al)$.
For each $d\in D$ the set $E_d = \{\al : d_\al = d\}$ is countable,
because every well-ordered subset of $\RR$ is countable.
Since $\xi = \bigcup_{d\in D} E_d$, we have
$|\xi| \le \max(|D|, \omega) = d(X)$.
\ep

The converse implication is not true in general.  For example,
if $X$ has the countable chain condition (ccc), then every
well-ordered chain in $C(X,\RR)$ is countable (because
$X \times \RR$ is also ccc).  However, the converse is true for
metric spaces:

\bl
If $(X, \varrho)$ is any non-empty metric space and $\prec$ is any total order
of the cardinal $d(X)$, then there is a chain in
$C(X,\RR)$ which is isomorphic to $\prec$.
\el
\bp
First, note that every countable total order is embeddable in $\RR$,
so if $d(X) = \omega$, then
the result follows trivially using constant functions.
In particular, we may assume that $X$ is infinite, and then
fix $D \subseteq X$ which is dense and of size $d(X)$.
For each $n\in\omega$, let $D_n$ be a subset of $D$ which is
maximal with respect to the property
$\forall d,e \in D_n \, [d \ne e \to \varrho(d,e) \ge 2^{2-n}]$.
Then $\bigcup_n D_n$ is also dense, so we may assume that
$\bigcup_n D_n = D$.
We may also assume that $\prec$ is a total order of the set $D$.
Now, we shall produce $f_d \in C(X,\RR)$ for $d\in D$ such
that $f_d < f_e$ whenever $d \prec e$.

For each $n$, if $c \in D_n$, define
$\varphi^n_c(x) = \max(0,\; 2^{-n} - \varrho(x,c))$.
For each $d \in D$, let
$\psi^n_d = \sum\{\varphi^n_c : c \in D_n  \ \&\ c \prec d\}$.
Since every $x \in X$ has a neighborhood on which all but
at most one of the $\varphi^n_c$ vanish, we have
$\psi^n_d \in C(X,[0,2^{-n}])$, and 
$\psi^n_d \le \psi^n_e$ whenever $d \prec e$.
Thus, if we let $f_d = \sum_{n<\omega}\psi^n_d$, we have
$f_d \in C(X,[0,2])$, and $f_d \le f_e$ whenever $d \prec e$.
But also, if $d \in D_n$ and $d \prec e$, then
$\psi^n_d(d) = 0 <  2^{-n} = \psi^n_e(d)$, so actually
$f_d < f_e$ whenever $d \prec e$.
\ep

Putting these lemmas together, we have:

\bt
\label{thm-woseqcont}
Let $(X,\varrho)$ be a metric space. Then there exists a well-ordered sequence
of length $\xi$ in $C(X,\RR)$ iff $\xi<d(X)^+$.
\et

\bcor
\label{cor-woseqcont}
A metric space $(X,\varrho)$ is separable iff every well-ordered sequence
in $C(X,\RR)$ is countable.
\ecor

\section{Sequences of Baire Class~1 Functions}
\label{sec-Baire}

If we replace continuous functions by Baire class 1 functions,
then Corollary \ref{cor-woseqcont} becomes false, since
on some separable metric spaces, we can get well-ordered sequences 
of every type less than $\omega_2$.
To prove this, we shall apply some basic facts
about $\subset^*$ on $\iP(\omega)$.
As usual, for $x,y \subseteq \omega$, we say that
$x \subseteq^* y$ iff $x \backslash y$ is finite.
Then $x \subset^* y$ iff $x \backslash y$ is finite 
and $y \backslash x$ is infinite.
This $\subset^*$ partially orders $\iP(\omega)$.

\bl
If $X \subset \iP(\omega)$ is a chain in the order $\subset^*$,
then on $X$ (viewed as a subset of the Cantor set
$2^\omega \cong \iP(\omega)$), there is a chain of Baire class 1 functions
which is isomorphic to $(X, \subset^*)$.
\el
\bp
Note that for each $x \in X$, 
$$
\{y \in X : y \subseteq^* x \} = \bigcup_{m\in \omega}
\{y \in X : \forall n \ge m \,[y(n) \le x(n)]\}\ \  ,
$$
which is an $F_\sigma$ set in $X$.  Likewise, the sets
$ \{y \in X :  y \supseteq^* x \}  $,
$ \{y \in X :  y \subset^* x  \} $, and
$\{y \in X :   y \supset^* x  \} $,
are all $F_\sigma$ sets in $X$,
and hence also $G_\delta$ sets.  It follows
that if $f_x : X \to \{0,1\}$ is the characteristic function
of $\{y\in X : y \subset^* x \}$, then $f_x : X \to \RR$
is a Baire class 1 function.
Then, $\{f_x : x \in X \}$ is the required chain.
\ep

\bl
For any infinite cardinal $\kappa$, suppose that $(\iP(\omega), \subset^*)$
contains a chain $\{x_\alpha : \alpha < \kappa\}$
(i.e., $\alpha < \beta \to  x_\alpha \subset^* x_\beta\}$).
Then $(\iP(\omega), \subset^*)$ contains a chain $X$ of size $\kappa$
such that every ordinal $\xi < \kappa^+$ is embeddable into $X$.
\el
\bp
Let $S = \bigcup_{1\le n < \omega} \kappa^n$.
For $s = (\alpha_1, \ldots , \alpha_{n-1}, \alpha_n) \in S$, let
$s^+ = (\alpha_1, \ldots , \alpha_{n-1}, \alpha_n+1)$.
Starting with the $x_{(\alpha)} = x_\alpha$,
choose $x_s \in \iP(\omega)$ by
induction on $\lh(s)$ so that
$x_s = x_{s\cat 0} \subset^*  x_{s\cat \alpha} \subset^*  x_{s\cat \beta}
\subset^* x_{s^+}$ whenever $s \in S$ and $0 < \alpha < \beta < \kappa$.
Let $X = \{x_s : s \in S\}$.  Then, whenever $x,y\in X$ with
$x \subset^* y$, the ordinal $\kappa$ is embeddable
in $(x,y) = \{z \in X : x \subset^* z \subset^* y\}$.
From this, one easily proves by induction on $\xi < \kappa^+$
(using $\cf(\xi)\leq\kappa$)
that $\xi$ is embeddable in each such interval $(x,y)$.
\ep

Since $\iP(\omega)$ certainly contains a chain of type $\omega_1$,
these two lemmas yield:

\bt
There is a separable metric space $X$ on which, for every
$\xi < \omega_2$, there is a well-ordered
chain of length $\xi$ of Baire class 1 functions.
\et

Under $\CH$, this is best possible, since there will be only
$2^\omega = \omega_1$ Baire class~1 functions on a \sms,
so there could not be a chain of length $\omega_2$.
Under $\neg \CH$, the existence of longer chains of Baire class 1 functions
depends on the model of set theory.  It is consistent with 
$\fc = 2^\omega$ being arbitrarily large that there is a chain
in $(\iP(\omega), \subset^*)$
of type $\fc$; for example, this is true under $\MA$ (see \cite{Do}).
In this case, there will be a separable $X$ with well-ordered chains
of all lengths less than $\fc^+$.
However, in the Cohen model, where $\fc$ can also be made
arbitrarily large, we never get chains of type $\omega_2$.
We shall prove this by using the following lemma, which relates it to
the rectangle problem:

\bl
\label{lemma-nochain}
Suppose that there is a separable metric space $Y$ 
with an $\omega_2$-chain of Borel subsets, $\{B_\alpha : \alpha < \omega_2\}$
(so, $\alpha < \beta \to B_\alpha \subsetneqq B_\beta$).
Then in $\omega_2 \times \omega_2$, the well-order relation $<$
is in the $\sigma$-algebra generated by the set of all rectangles,
$\{ S\times T : S,T \in \iP(\omega_2) \}$.
\el
\bp
Each $B_\alpha$ has some countable Borel rank.  Since there are only
$\omega_1$ ranks, we may, by passing to a subsequence, assume
that the ranks are bounded.  Say, each $B_\alpha$ is a $\Sigma^0_\mu$
set for some fixed $\mu < \omega_1$.  

Let $J = \omega^\omega$, and let $A \subseteq Y\times J$ be\
a universal $\Sigma^0_\mu$ set; that is, $A$ is $\Sigma^0_\mu$
in $Y \times J$ and every $\Sigma^0_\mu$ subset of $Y$
is of the form $A^j = \{y : (y,j) \in A\}$ for some $j \in J$
(see \cite{Ku}, \S31).
Now, for $\alpha,\beta < \omega_2$,
fix $y_\alpha \in B_{\alpha+1} \backslash B_\alpha$, and
fix $j_\beta \in J$ such that
$A^{j_\beta} = B_{\beta}$.
Then $\alpha < \beta$ iff $(y_\alpha, j_\beta) \in A$.
Thus, $\{(y_\alpha, j_\beta): \alpha < \beta < \omega_2\}$
is a Borel subset of
$\{y_\alpha : \alpha < \omega_2\} \times \{j_\beta : \beta < \omega_2\}$,
and is hence in the $\sigma$-algebra generated by open rectangles,
so $<$, as a subset of $\omega_2\times\omega_2$,
is in the $\sigma$-algebra generated by rectangles.
\ep

\bt
Assume that the well-order relation $<$ on $\omega_2$
is not in the $\sigma$-algebra generated by the set of all rectangles.
Then no separable metric space can have a 
chain of length $\omega_2$ of Baire class 1 functions.
\et
\bp
Suppose that $\{f_\al : \al<\omega_2 \}$ is a chain 
of Baire class one functions on the \sms\ $X$.
Let $B_\alpha = \{(x,r) \in X\times\RR : r \le f_\alpha(x)\}$.
Then the $B_\alpha$ form an
$\omega_2$-chain of Borel subsets of the \sms\  $X\times \RR$,
so we have a contradiction by Lemma \ref{lemma-nochain}.
\ep

Finally, we point out that the hypothesis of this theorem
is consistent, since it holds in the extension
$V[G]$ formed by adding $\ge \omega_2$ Cohen reals to a ground
model $V$ which satisfies $\CH$. 
This fact was first proved in \cite{Kun}.
It also follows from the more general principle $HP_2(\omega_2)$
of  Brendle,  Fuchino, and  Soukup \cite{BFS}.
They define this principle, prove that it holds in
Cohen extensions (and in a number of other forcing extensions),
and show the following:

\bl
\label{lemma-hp}
$HP_2(\kappa)$ implies that if $R$ is any relation on
$\iP(\omega)$ which is first-order
definable over $H(\omega_1)$ from a
fixed element of $H(\omega_1)$, then there is no $X \subseteq \iP(\omega)$
such that $(X; R)$ is isomorphic to $(\kappa; <)$.
\el

These matters are also
discussed in \cite{JK}, which indicates how
such statements are verified in Cohen extensions.
Here, $H(\omega_1)$ denotes the set of hereditarily countable
sets.

\bl
$HP_2(\omega_2)$ implies that 
in $\omega_2 \times \omega_2$, the well-order relation $<$
is not in the $\sigma$-algebra generated by the set of all rectangles,
$\{ S\times T : S,T \in \iP(\omega_2) \}$.
\el
\bp
Suppose that  $<$ were in this $\sigma$-algebra.
Then we would have fixed $K_n \subseteq \omega_2$ for $n < \omega$
such that $<$ is in the $\sigma$-algebra generated by all the
$K_m \times K_n$.  For each $\alpha$, let
$u_\alpha = \{n\in\omega : \alpha\in K_n\}$.
There is then a formula $\varphi(x,y,z)$ and a fixed
$w \in H(\omega_1)$ such that 
for all $\alpha, \beta < \omega_2$,
$\; \alpha < \beta$ iff $H(\omega_1) \models \varphi(u_\alpha, u_\beta, w)$;
here, $w$ encodes the
particular countable boolean combination used to get $<$ from
the $K_n$.
Now, if $X = \{u_\alpha : \alpha < \omega_2\}$, then
$\varphi$ defines a relation $R$ on $H(\omega_1)$ such that
$(X; R)$ is isomorphic to $(\omega_2; <)$, contradicting
Lemma \ref{lemma-hp}.
\ep

\bigskip

\noindent
\textsc{Department of Analysis, E\"otv\"os Lor\'and University,
Budapest, P\'azm\'any P\'eter s\'et\'any 1/c, 1117, Hungary}

\textit{Email address}: \verb+emarci@cs.elte.hu+

\bigskip

\noindent
\textsc{Department of Mathematics,
University of Wisconsin, Madison, WI 53706, USA}

\textit{Email address}: \verb+kunen@math.wisc.edu+

\textit{URL:} \verb+http://www.math.wisc.edu/~kunen+


\begin{thebibliography}{99}

\bibitem{BFS}
J. Brendle, S. Fuchino, and L. Soukup,
Coloring ordinals by reals, 
to appear.


\bibitem{Do} E. K. van Douwen, 
The integers and topology, in
\textsl{Handbook of Set-Theoretic Topology},
North-Holland, Amsterdam, 1984, pp. 111-167. 

\bibitem{El} M. Elekes, Linearly ordered families of
Baire~1 functions, \textsl{Real Analysis Exchange}, to appear.

\bibitem{JK} I. Juh\'asz and K. Kunen,
The power set of $\omega$,
\textsl{Fundamenta Mathematicae}, Vol \textbf{170} (2001), 257-265.

\bibitem{Ko} P. Komj\'ath, Ordered families of
Baire-2-functions, \textsl{Real Analysis Exchange}, Vol \textbf{15}
(1989-90), 442-444.

\bibitem{Kun} K. Kunen, \textsl{Inaccessibility Properties of Cardinals},
Doctoral Dissertation, Stanford, 1968.

\bibitem{Ku} K. Kuratowski, \textsl{Topology, Vol. 1}, Academic Press, 1966.

\end{thebibliography}
\end{document}